\documentclass[11pt]{amsart}

\usepackage{latexsym,rawfonts}
\usepackage{amsfonts,amssymb}
\usepackage{amsmath,amsthm}
\usepackage{mathrsfs}

\usepackage[plainpages=false]{hyperref}
\usepackage{graphicx}
\usepackage{color}

\newtheorem{theorem}{Theorem}[section]
\newtheorem{corollary}[theorem]{Corollary}
\newtheorem{lemma}[theorem]{Lemma}

\theoremstyle{definition}

\numberwithin{equation}{section}

\DeclareMathOperator{\DIV}{div}
\DeclareMathOperator{\VOL}{vol}

\DeclareMathOperator{\TR}{tr}

\DeclareMathOperator{\DIAG}{diag}

\DeclareMathOperator{\SPAN}{span}

\newcommand{\R}{\mathbb{R}}

\newcommand{\dd}{\mathop{}\!\mathrm{d}}

\newcommand{\norm}[1]{\left\Vert#1\right\Vert}

\newcommand{\one}{\mathbf{1}}

\newcommand{\uS}{\mathbb{S}^{n-1}}

\begin{document}

\title{Remarks on uniqueness and nonuniqueness for the isotropic $L_p$-Minkowski problem with $p<-n$}

\author{YanNan Liu}
\address{School of Mathematics and Statistics,
  Beijing Technology and Business University, 
  Beijing 100048, P.R. China}
\email{liuyn@th.btbu.edu.cn}

\author{Jian Lu}
\address{School of Mathematics and Statistics, Nanjing University of Science and Technology, 
  Nanjing 210094, P.R. China}
\email{lj-tshu04@163.com}

\thanks{The author Liu was supported by Natural Science Foundation of China (12141103).
The author Lu was supported by
Natural Science Foundation of China (12671241, 12526435), and by
Fundamental Research Funds for Central Universities (30926WZCT202601009).
}

\date{}

\subjclass[2020]{52A20, 35J96, 35A02}

\keywords{$L_p$-Minkowski problem, Monge-Amp\`ere equation, uniqueness, nonuniqueness}

\begin{abstract}
We study uniqueness and nonuniqueness for the isotropic $L_p$-Minkowski problem with $p<-n$.
For $n\geq3$, we prove uniqueness for $p$ just below $-n$ without any symmetry assumptions.
Moreover, for $n\geq8$, we obtain nonuniqueness for $p$ just above $-2n-3$, even under full simplex symmetry.
The proofs use the Chou-Wang obstruction and Lyapunov-Schmidt reduction.
\end{abstract}

\maketitle

\section{Introduction}

Let $\uS$ be the unit sphere in the $n$-dimensional Euclidean space $\R^n$.
We consider 
the following isotropic $L_p$-Minkowski problem
\begin{equation}\label{LpMP1}
 \det(\nabla^2h+hI)=h^{p-1}\qquad\text{on }\uS,
\end{equation}
where $h$ is the support function of a convex body $K\subset\R^n$, $\nabla^2h$ is the covariant Hessian matrix in a local orthonormal frame,
$I$ is the $(n-1)\times(n-1)$ identity matrix, 
and $p$ is a real number.

The $L_p$-Minkowski problem, introduced by Lutwak \cite{Lutwak}, asks for a convex body with prescribed $L_p$ surface area measure $h^{1-p}\dd S_K$.
Here $\dd S_K$ is the surface area measure of $K$.
It includes the classical Minkowski problem when $p=1$, and \eqref{LpMP1} is the case of constant prescribed density.
For existence results, we refer to \cite{BBCY,BLYZ,CW,GLW,HLYZ,LW,Lutwak,LYZ}, and for regularity to \cite{BBC,CW,LO}.
For uniqueness results with constant or nonconstant prescribed data, see \cite{BS,BCD,CHLL,HLX}; for nonuniqueness, see \cite{BS,CLZ,DLL,JLW}.

The constant function $h\equiv1$ solves \eqref{LpMP1} for every $p$.
In this paper,
we are mainly concerned with whether it is the only positive solution when $n\geq3$.
Some basic uniqueness and nonuniqueness results are summarized in Table~\ref{tab:1}.
\begin{table}
\caption{Uniqueness/nonuniqueness of Eq. \eqref{LpMP1} for $n\geq3$}\label{tab:1}
\begin{tabular}{@{}lll@{}}
\hline
Exponent & Positive solutions & Reference \\
\hline
$p>1$ & $h\equiv1$ (up to dilation if $p=n$) & Lutwak \cite{Lutwak} \\
$p=1$ & $h\equiv1$ (up to translation) & Lutwak-Yang-Zhang \cite{LYZ} \\
$-n<p<1$ & $h\equiv1$ & Brendle-Choi-Daskalopoulos \cite{BCD} \\
$p=-n$ & Ellipsoids; see \eqref{eq:13} & Calabi \cite{Calabi} \\
$p<-2n-3$ & Nonconstant solutions exist & Du-Liu-Lu \cite{DLL} \\
\hline
\end{tabular}
\end{table}

The uniqueness question remains open for $p\in[-2n-3,-n)$.
We give partial answers near the two endpoints of this interval.

Our first result establishes uniqueness for $p$ sufficiently close to $-n$ from below.

\begin{theorem}\label{thm:1}
Let $n\geq3$.
There exists $\varepsilon_n\in(0,1]$ such that, for every
$p\in(-n-\varepsilon_n,-n)$, Eq. \eqref{LpMP1} has a unique
positive solution, namely $h\equiv1$.
\end{theorem}

Our second result establishes nonuniqueness for $p$ sufficiently close to $-2n-3$ from above, even under the full symmetry of a regular simplex.
\begin{theorem}\label{thm:2}
Let $n\ge8$. 
There exists $\beta_n\in(0,n+3)$ such that,
for every 
\[ p\in(-2n-3,-2n-3+\beta_n), \]
Eq. \eqref{LpMP1} has a nonconstant positive solution,
which is invariant under the symmetry group of a regular $n$-simplex centered at the origin.
\end{theorem}

To the best of our knowledge, Theorem~\ref{thm:1} is the first uniqueness result for the isotropic $L_p$-Minkowski problem with $p<-n$ in dimensions $n\geq3$ without any symmetry assumptions.
While Theorem~\ref{thm:2} shows that nonuniqueness can occur even under highly strong symmetry restrictions, illustrating the complexity of the uniqueness question for $p<-n$.

The paper is organized as follows.
In Section~\ref{sec:2}, we recall some basic facts used in the proofs.
In Section~\ref{sec:3}, we prove Theorem~\ref{thm:1} through a careful application of the Chou-Wang obstruction; see Corollary \ref{cor:CW-1}.
Finally, in Section~\ref{sec:4}, we apply Lyapunov-Schmidt reduction to prove Theorem~\ref{thm:2}.

\section{Preliminaries}\label{sec:2}

In this section, we recall some facts needed in the following proofs.
We write
\[
 A_h=\nabla^2h+hI
\]
in a local orthonormal frame, as in \eqref{LpMP1}.
Every positive support function $h$ solving \eqref{LpMP1} is smooth and satisfies $A_h>0$; see \cite[Proposition~1.2 and the remark following it]{CW}.

If $K$ has a positive smooth support function $h$ with $A_h>0$, its inverse Gauss parametrization is
\begin{equation}\label{eq:2}
 X(x)=h(x)x+\nabla h(x),\qquad
 \dd S_K=\det A_h\,\dd x.
\end{equation}
In fact, $\dd X$ restricted to the tangent space is $A_h$, which gives the second formula.
The polar body
\[
 K^\circ=\{z\in\R^n:z\cdot y\leq1\text{ for all }y\in K\}
\]
has radial function $\rho_{K^\circ}(x)=h(x)^{-1}$.

For solution $h$ to Eq. \eqref{LpMP1}, there are two useful identities.
The first one is
\begin{equation}\label{eq:3}
 \int_{\uS}x h^{p-1}\,\dd x=0,
\end{equation}
which directly follows the divergence theorem.
The second one is a special case of the Chou-Wang obstruction \cite{CW}.

\begin{lemma} \label{lem:CW-f}
Assume $p\neq0$.
Let $h\in C^3(\uS)$ be a solution to the $L_p$ Minkowski problem
\begin{equation}\label{LpMPf}
 \det(\nabla^2h+hI)=fh^{p-1}\qquad\text{on }\uS.
\end{equation}
Then, for every projective vector field $\xi$, we have
\begin{equation}\label{CW-f}
\int_{\mathbb{S}^{n-1}} h^p \left(\nabla_\xi f + \frac{p+n}{n} f\DIV\xi\right) 
 = 0.
\end{equation}
Here $\xi$ is given by
\[
\xi(x) \equiv Bx - (x^T Bx)x,
 \qquad x \in \mathbb{S}^{n-1},
\]
where $B$ is any $n\times n$ matrix.
\end{lemma}

Note that a generalization of \eqref{CW-f} has been obtained in \cite{DLL}.

\begin{proof}[Proof of Lemma \ref{lem:CW-f}]
We rewrite \eqref{LpMPf} as
\[
 h^{1+n}\det(\nabla^2h+hI)=fh^{p+n}=:\widetilde f.
\]
Applying the Chou-Wang obstruction \cite[Proposition~F]{CW} for exponent $-n$, we have
\begin{equation*}
\begin{split}
 0
 =\int_{\uS}h^{-n}\nabla_\xi\widetilde f
 &=\int_{\uS}h^p\nabla_\xi f
       +\frac{p+n}{p}\int_{\uS}f\nabla_\xi(h^p)\\
 &=\int_{\uS}h^p\nabla_\xi f
       -\frac{p+n}{p}\int_{\uS}h^p
       \left(\nabla_\xi f+f\DIV\xi\right)\\
 &=-\frac{n}{p}\int_{\uS}h^p
       \left(\nabla_\xi f+\frac{p+n}{n}f\DIV\xi\right),
\end{split}
\end{equation*}
which proves \eqref{CW-f}.
\end{proof}

For the isotropic case, the Chou-Wang obstruction becomes very simple.

\begin{corollary}\label{cor:CW-1}
Assume $p\notin\{0,-n\}$.
Let $h\in C^3(\uS)$ solve Eq. \eqref{LpMP1}.
Then we have
\begin{equation}\label{CW-1}
 \int_{\uS}h^p x_i x_j
 =\frac{\delta_{ij}}{n}\int_{\uS}h^p,
 \qquad \forall\,1\le i,j\le n.
\end{equation}
\end{corollary}

\begin{proof}
By Lemma~\ref{lem:CW-f} with $f\equiv1$, we have
\[
 \frac{p+n}{n}\int_{\uS}h^p\DIV\xi=0.
\]
Since $p\ne-n$ and $\DIV\xi=\TR B-nx^TBx$, we obtain
\[
 n\int_{\uS}h^p \,x^TBx
=\TR B \int_{\uS}h^p.
\]
Since $B$ is arbitrary, this proves \eqref{CW-1}.
\end{proof}

At the critical exponent $p=-n$, we shall use the following classification of Calabi \cite{Calabi}.

\begin{lemma}\label{lem:5}
The support functions satisfying
\[
 \det(\nabla^2h+hI)=h^{-n-1}\qquad\text{on }\uS
\]
are given by
\begin{equation}\label{eq:13}
 h(x)=(x^TQx)^{1/2},\qquad \det Q=1,
\end{equation}
where $Q$ is symmetric and positive definite.
\end{lemma}

Finally, for the linearized equations, we use the convention $\Delta=\TR\nabla^2$ and recall that
\begin{equation}\label{eq:1}
 \Delta Y=-\ell(\ell+n-2)Y,\qquad Y\in\mathcal H_\ell,
\end{equation}
where $\mathcal H_\ell$ is the space of spherical harmonics of degree $\ell$.

\section{Uniqueness}\label{sec:3}

Throughout this section, $n\geq3$.
We first state the two assertions needed for the proof.

\begin{lemma}\label{lem:3}
Let $p_j<-n$ tend to $-n$, and let $h_j$ be positive solutions of \eqref{LpMP1} with exponents $p_j$.
Then $h_j\to1$ in $C^k(\uS)$ for every integer $k\geq0$.
\end{lemma}

\begin{lemma}\label{lem:1}
Fix $\alpha\in(0,1)$.
There exists $\rho=\rho(n,\alpha)>0$ such that, if $p\in[-n-1,-n)$ and a positive solution $h$ of \eqref{LpMP1} satisfies
\[
 \norm{h-1}_{C^{2,\alpha}(\uS)}<\rho,
\]
then $h\equiv1$.
\end{lemma}

Assuming these for the moment, we prove the theorem.

\begin{proof}[Proof of Theorem~\ref{thm:1}]
The function $h\equiv1$ solves \eqref{LpMP1}.
If the asserted uniqueness failed, for every integer $j\geq1$ there would be an exponent
$p_j\in(-n-1/j,-n)$ and a positive solution $h_j\not\equiv1$.
By Lemma~\ref{lem:3}, $h_j\to1$ in $C^3(\uS)$, and hence in $C^{2,\alpha}(\uS)$ for any fixed $\alpha\in(0,1)$.
Lemma~\ref{lem:1} then gives $h_j\equiv1$ for all sufficiently large $j$, a contradiction.
Thus the required neighborhood exists, and we may decrease its size so that $\varepsilon_n\leq1$.
\end{proof}

It remains to prove Lemma~\ref{lem:3} and Lemma~\ref{lem:1}.

\subsection{Uniform estimates}
We begin with an estimate for convex bodies satisfying the two moment conditions.

\begin{lemma}\label{lem:2}
Let $L\subset\R^n$ be a convex body containing the origin in its interior, and let $\delta\in[0,1]$.
Set
\[
 Z_L=\int_L|y|^\delta\,\dd y,\qquad
 \dd\mu(y)=Z_L^{-1}|y|^\delta\,\dd y.
\]
Suppose that
\begin{equation}\label{eq:4}
 \int_L y\,\dd\mu=0,\qquad
 \int_L\frac{y y^T}{|y|^2}\,\dd\mu=\frac1nI_n.
\end{equation}
If $R=\max_{y\in L}|y|$, then
\begin{equation}\label{eq:5}
 c_*RB^n\subset L\subset RB^n,\qquad c_*=\frac{1}{32n^3 4^n},
\end{equation}
where $B^n$ is the Euclidean unit ball.
The value of the second integrand at $y=0$ is immaterial.
\end{lemma}

\begin{proof}
The assumptions are preserved by dilation, so we may assume $R=1$.
Choose $z\in L$ with $|z|=1$.
By convexity, $\frac34z+\frac14L\subset L$, and every point of this set has norm at least $1/2$.
Since $0\leq\delta\leq1$, we have
\begin{equation}\label{eq:6}
 Z_L\geq\frac12\,4^{-n}\VOL(L).
\end{equation}
Let $[a,b]$ be the projection of $L$ in the direction $z$, and let $m(t)$ be the $(n-1)$-dimensional volume of $L\cap\{y:y\cdot z=t\}$.
We have $b=1$ and $a<0$, so $b-a\geq1$.
For $a<t<b$, join this section to points of $L$ at heights $a$ and $b$.
The two resulting cones lie in $L$, their interiors are disjoint, and their volumes sum to $(b-a)m(t)/n$.
Consequently,
\[
 m(t)\leq\frac{n}{b-a}\VOL(L)\leq n\VOL(L).
\]
Integration over $-r\leq t\leq r$ gives $\VOL(L\cap rB^n)\leq2nr\VOL(L)$.
Combining this with \eqref{eq:6} and $|y|^\delta\leq1$ on $L$, we obtain
\[
 \mu(L\cap rB^n)\leq C_*r,\qquad C_*=4^{n+1}n.
\]
Put $r_*=(2nC_*)^{-1}$.
For any unit vector $x$, the second identity in \eqref{eq:4} gives
\[
 \begin{split}
 \frac1n
 &=\int_L\frac{(x\cdot y)^2}{|y|^2}\,\dd\mu\\
 &\leq\mu(L\cap r_*B^n)+\frac1{r_*}\int_L|x\cdot y|\,\dd\mu \\
 &\leq\frac1{2n}+\frac1{r_*}\int_L|x\cdot y|\,\dd\mu.
 \end{split}
\]
Here we used $(x\cdot y)^2/|y|^2\leq|x\cdot y|/r_*$ when $|y|\geq r_*$.
The first identity in \eqref{eq:4} now yields
\[
 h_L(x)\geq\int_L(x\cdot y)_+\,\dd\mu
 =\frac12\int_L|x\cdot y|\,\dd\mu
 \geq\frac{r_*}{4n}=c_*.
\]
Thus every supporting half-space of $L$ contains $c_*B^n$, which proves \eqref{eq:5}.
\end{proof}

We next apply this lemma to the polar body and then estimate the principal radii of curvature.

\begin{lemma}\label{lem:4}
There are constants $\lambda,\Lambda>0$, depending only on $n$, such that every positive solution of \eqref{LpMP1} with $p\in[-n-1,-n)$ satisfies
\begin{equation}\label{eq:7}
 c_*\leq h\leq c_*^{-1},\qquad h^2+|\nabla h|^2\leq c_*^{-2},
\end{equation}
and
\begin{equation}\label{eq:8}
 \lambda I\leq A_h\leq\Lambda I.
\end{equation}
Moreover, for every integer $k\geq0$, there is $C_k=C_k(n)$ such that $\norm{h}_{C^k(\uS)}\leq C_k$.
Here $c_*$ is as in Lemma~\ref{lem:2}.
\end{lemma}

\begin{proof}
Write $p=-n-\delta$, where $0<\delta\leq1$.
Let $K$ be the convex body with support function $h$, set $L=K^\circ$, and write $R=\max_{y\in L}|y|$.
By polar coordinates and $\rho_L=h^{-1}$,
\[
 \begin{split}
 Z_L&=\frac1{n+\delta}\int_{\uS}h^{p} \dd x,\\
 \int_L y|y|^\delta \dd y
 &=\frac1{n+\delta+1}\int_{\uS}x h^{p-1} \dd x,\\
 \int_L\frac{y y^T}{|y|^2}|y|^\delta \dd y
 &=\frac1{n+\delta}\int_{\uS}x x^T h^{p} \dd x.
 \end{split}
\]
Thus \eqref{eq:3} and Corollary~\ref{cor:CW-1} give \eqref{eq:4}.
By Lemma~\ref{lem:2} and polarity,
\[
 R^{-1}B^n\subset K\subset(c_*R)^{-1}B^n.
\]
In particular, if $m=\min h$ and $M=\max h$, then $M/m\leq c_*^{-1}$.
At points where $h$ attains its minimum and maximum, respectively, \eqref{LpMP1} gives
\[
 m^{p-1}\geq m^{n-1},\qquad M^{p-1}\leq M^{n-1}.
\]
Since $p<n$, we have $m\leq1\leq M$.
It follows that $c_*\leq m\leq M\leq c_*^{-1}$.
Moreover, $K\subset c_*^{-1}B^n$, and \eqref{eq:2} gives
$h^2+|\nabla h|^2=|X|^2\leq c_*^{-2}$.
This proves \eqref{eq:7}.

For the second derivative estimate \eqref{eq:8}, we first bound the trace of $A_h$ from above.
Since $A_h>0$, this bounds every eigenvalue from above; the equation and \eqref{eq:7} then give a lower bound for their product and hence for each eigenvalue.
For simplification, write
\[
 A_{ij}=h_{,ij}+h\delta_{ij},
 \qquad (A^{ij})=(A_{ij})^{-1},
 \qquad \mathcal L=A^{ij}\nabla_i\nabla_j.
\]
Here repeated tangent indices are summed.

We first compute $\mathcal L(\TR A)$.
Note that $\TR A=\Delta h+(n-1)h$.
The tensor $A$ is Codazzi, that is, $\nabla_k A_{ij}=\nabla_i A_{kj}$.
Commuting these derivatives on the unit sphere and contracting gives
\begin{equation}\label{eq:9}
 \Delta A_{ij}=(\TR A)_{,ij}+(n-1)A_{ij}-(\TR A)\delta_{ij}.
\end{equation}
To evaluate the contraction with $A^{-1}$, we differentiate $\log\det A=(p-1)\log h$ twice and sum over the derivative directions to obtain
\[
 A^{ij}\Delta A_{ij}=(p-1)\Delta\log h+E,\qquad
 E=\sum_k\TR\big[(A^{-1/2}(\nabla_k A)A^{-1/2})^2\big]\geq0.
\]
By \eqref{eq:9} and $A^{ij}A_{ij}=n-1$, we obtain
\begin{equation}\label{eq:10}
 \begin{split}
 \mathcal L(\TR A)
 &=(p-1)\Delta\log h+E-(n-1)^2+(\TR A)\TR(A^{-1})\\
 &=(p-1)\left(\frac{\TR A}{h}-(n-1)-\frac{|\nabla h|^2}{h^2}\right)
       +E-(n-1)^2+(\TR A)\TR(A^{-1})\\
 &\geq (\TR A)\TR(A^{-1})-C_1\TR A-(n-1)^2,\qquad C_1=\frac{n+2}{c_*}.
 \end{split}
\end{equation}
We used $E\geq0$, $p-1<0$, and \eqref{eq:7} in the last line.

To compensate for the negative term $-C_1\TR A$ in \eqref{eq:10}, we shall add a suitable multiple of $|X|^2$.
We have $(|X|^2)_{,i}=2A_{ij}h_{,j}$.
The Codazzi identity and $A^{ij}\nabla_k A_{ij}=(p-1)h_{,k}/h$ therefore give
\begin{equation}\label{eq:11}
 \begin{split}
 \mathcal L(|X|^2)
 &=2A^{ij}(\nabla_j A_{ik})h_{,k}
       +2A^{ij}A_{ik}h_{,kj}\\
 &=2(p-1)\frac{|\nabla h|^2}{h}+2\Delta h \\
 &=2\TR A-2(n-1)h+2(p-1)\frac{|\nabla h|^2}{h}\\
 &\geq2\TR A-C_2,\qquad
 C_2=\frac{2(n-1)}{c_*}+\frac{2(n+2)}{c_*^3}.
 \end{split}
\end{equation}

Now, let $x_0$ be a maximum point of $\TR A+\frac{C_1+1}{2}|X|^2$.
Since $A>0$, the maximum principle and \eqref{eq:10}--\eqref{eq:11} give, at $x_0$,
\[
 \begin{split}
 0&\geq\mathcal L\left(\TR A+\frac{C_1+1}{2}|X|^2\right)\\
 &\geq (\TR A)\TR(A^{-1})+\TR A-(n-1)^2-\frac{C_1+1}{2}C_2.
 \end{split}
\]
In particular, $\TR A(x_0)\leq (n-1)^2+\frac{C_1+1}{2}C_2$.
Since $0\leq |X|^2\leq c_*^{-2}$ by \eqref{eq:7}, the definition of $x_0$ gives, for every $x\in\uS$,
\begin{equation}\label{eq:12}
 \begin{split}
 \TR A(x)&\leq \TR A(x)+\frac{C_1+1}{2}|X(x)|^2\\
 &\leq \TR A(x_0)+\frac{C_1+1}{2}|X(x_0)|^2\\
 &\leq (n-1)^2+\frac{C_1+1}{2}C_2+\frac{C_1+1}{2}c_*^{-2}=:C_3.
 \end{split}
\end{equation}

Each eigenvalue of $A$ is positive and at most $\TR A\leq C_3$.
To bound the smallest eigenvalue from below, we use the determinant.
By \eqref{eq:7} and the range of $p$, we have $\det A=h^{p-1}\geq c_*^{n+2}$.
Dividing by the product of the other $n-2$ eigenvalues gives a lower bound of $c_*^{n+2}C_3^{-(n-2)}$ for the smallest eigenvalue.
This proves \eqref{eq:8}.

Finally, \eqref{eq:7}--\eqref{eq:8} and the standard Evans--Krylov and Schauder estimates give uniform bounds in every $C^k$ norm; see \cite[Chapters~6 and~17]{GT}.
The constants depend only on $n$ and $k$, since $p$ lies in a bounded interval and all preceding bounds depend only on $n$.
\end{proof}

\subsection{The limit at the critical exponent}

\begin{proof}[Proof of Lemma~\ref{lem:3}]
After discarding finitely many terms, we have $p_j\in[-n-1,-n)$.
By Lemma~\ref{lem:4} and a diagonal compactness argument, every subsequence has a further subsequence converging in every $C^k$ norm to a smooth function $h_\infty$.
The bounds \eqref{eq:7}--\eqref{eq:8} persist in the limit, so $h_\infty>0$, $A_{h_\infty}>0$, and passage to the limit in \eqref{LpMP1} gives
\[
 \det A_{h_\infty}=h_\infty^{-n-1}.
\]
By Lemma~\ref{lem:5},
\[
 h_\infty(x)=(x^TQx)^{1/2},\qquad Q>0,\qquad \det Q=1.
\]
Moreover, $h_j^{p_j}\to h_\infty^{-n}$ uniformly, since $h_j$ are uniformly bounded away from zero.
Corollary~\ref{cor:CW-1} consequently gives
\begin{equation}\label{eq:14}
 \int_{\uS}h_\infty^{-n}x_i x_j
 =\frac{\delta_{ij}}{n}\int_{\uS}h_\infty^{-n},
 \qquad \forall\,1\le i,j\le n.
\end{equation}

We show that this identity forces $Q=I_n$.
After an orthogonal change of coordinates, we have
\[
 Q=\DIAG(\lambda_1,\ldots,\lambda_n),\qquad h_\infty(x)^2=x^TQx=\sum_k\lambda_kx_k^2.
\]
If $\lambda_i>\lambda_j$, let $\tau$ interchange coordinates $i$ and $j$, and put $\varphi(x)=x_i^2-x_j^2$.
Then
\begin{equation}\label{eq:19}
 h_\infty^2(x)-h_\infty^2(\tau x)=(\lambda_i-\lambda_j)\varphi(x),\qquad
 \varphi(\tau x)=-\varphi(x).
\end{equation}
Since $\tau$ preserves spherical measure, the change of variables $x=\tau y$ gives
\[
 \begin{split}
 \int_{\uS}h_\infty^{-n}(x)\varphi(x)\,\dd x
 &=\int_{\uS}h_\infty^{-n}(\tau y)\varphi(\tau y)\,\dd y\\
 &=-\int_{\uS}h_\infty^{-n}(\tau x)\varphi(x)\,\dd x.
 \end{split}
\]
Therefore, we have
\begin{equation}\label{eq:20}
 2\int_{\uS}h_\infty^{-n}(x) \varphi(x) \dd x
 =\int_{\uS}\big(h_\infty^{-n}(x)-h_\infty^{-n}(\tau x)\big) \varphi(x) \dd x<0.
\end{equation}
Indeed, note that 
$h_\infty(x)-h_\infty(\tau x)$ has the same sign as $\varphi(x)$ due to \eqref{eq:19}.
Since $t\mapsto t^{-n}$ is strictly decreasing on $(0,\infty)$, the difference $h_\infty^{-n}(x)-h_\infty^{-n}(\tau x)$ has the opposite sign whenever $\varphi(x)\ne0$.
Thus the integrand is nonpositive everywhere and strictly negative on $\{x\in\uS:\varphi(x)\ne0\}$.
Thus, \eqref{eq:20} holds.

On the other hand, by \eqref{eq:14}, there is
\[
 \int_{\uS}h_\infty^{-n}(x)\varphi(x)\,\dd x
 =\int_{\uS}h_\infty^{-n}(x)(x_i^2-x_j^2)\,\dd x=0,
\]
a contradiction with \eqref{eq:20}.
All eigenvalues of $Q$ are therefore equal, and $\det Q=1$ gives $Q=I_n$.
Thus every subsequential limit is $h_\infty=1$.
If convergence to $1$ failed in some $C^k$ norm, a subsequence staying a fixed positive distance from $1$ in that norm would have a further subsequence converging smoothly to $1$, a contradiction.
This proves the lemma.
\end{proof}

\subsection{Local rigidity}

The linearization at $h=1$, $p=-n$ has kernel $\mathcal H_2$.
We use \eqref{CW-1} to estimate this component before estimating its orthogonal complement.

\begin{proof}[Proof of Lemma~\ref{lem:1}]
Let $v=h-1$, and let $\Pi_2$ be the $L^2(\dd x)$-orthogonal projection onto $\mathcal H_2$.
This space consists of $x\mapsto x^TTx$, where $T$ is symmetric and trace-free.
Hence \eqref{CW-1} is equivalent to
$\Pi_2((1+v)^p)=0$.
Since $\Pi_2(1)=0$, we have
\[
 p\Pi_2v=-\Pi_2\big((1+v)^p-1-pv\big).
\]
For $\norm{v}_{C^0}\leq1/2$, Taylor's formula gives
\[
 (1+v)^p-1-pv
 =p(p-1)v^2\int_0^1(1-t)(1+tv)^{p-2}\,\dd t.
\]
The coefficient on the right is bounded uniformly for $p\in[-n-1,-n]$.
The finite-rank projection $\Pi_2$ is bounded from $C^0$ to $C^{2,\alpha}$: this follows by expanding in a fixed orthonormal basis of the smooth functions in $\mathcal H_2$.
Using $|p|\geq n$, we conclude that
\begin{equation}\label{eq:15}
 \norm{\Pi_2v}_{C^{2,\alpha}}
 \leq C\norm{v}_{C^0}^2
 \leq C\norm{v}_{C^{2,\alpha}}^2.
\end{equation}
Here and below $C$ depends only on $n$ and $\alpha$.

For $\norm{v}_{C^{2,\alpha}}$ sufficiently small, the logarithmic form of \eqref{LpMP1} is
\begin{equation}\label{eq:16}
 0=\log\det(I+A_v)+(1-p)\log(1+v)=L_pv+N_p(v),
\end{equation}
where $L_p=\Delta+n-p$, and
\[
 N_p(v)
 =\log\det(I+A_v)-\TR A_v
  +(1-p)\big(\log(1+v)-v\big).
\]
Indeed, $\TR A_v+(1-p)v=\Delta v+(n-p)v=L_pv$, so $N_p(v)$ is the remainder after subtracting the linear terms.
We claim that
\begin{equation}\label{eq:17}
 \norm{N_p(v)}_{C^{0,\alpha}}\leq C\norm{v}_{C^{2,\alpha}}^2
\end{equation}
uniformly for $p\in[-n-1,-n]$.
To verify this, Taylor's formula with integral remainder gives
\[
 \begin{split}
 \log\det(I+A_v)-\TR A_v
 &=-\int_0^1(1-t)\TR\big[((I+tA_v)^{-1}A_v)^2\big]\,\dd t,\\
 \log(1+v)-v&=-v^2\int_0^1\frac{1-t}{(1+tv)^2}\,\dd t.
 \end{split}
\]
On a sufficiently small $C^{2,\alpha}$ ball, we have $\norm{v}_{C^0}\leq1/2$ and $\norm{A_v}_{C^0}\leq1/2$.
Thus $1+tv\geq1/2$ and $I+tA_v\geq\frac12I$ for $0\leq t\leq1$.
Smoothness of inversion on this range gives
\[
 \sup_{0\leq t\leq1}\left(
 \norm{(I+tA_v)^{-1}}_{C^{0,\alpha}}
 +\norm{(1+tv)^{-1}}_{C^{0,\alpha}}\right)\leq C.
\]
Each integral remainder therefore contains two factors of $A_v$ or $v$, with all other factors uniformly bounded in $C^{0,\alpha}$.
Since $|1-p|\leq n+2$, the H\"older product estimate yields
\[
 \norm{N_p(v)}_{C^{0,\alpha}}
 \leq C\left(\norm{A_v}_{C^{0,\alpha}}^2+\norm{v}_{C^{0,\alpha}}^2\right)
 \leq C\norm{v}_{C^{2,\alpha}}^2,
\]
which proves \eqref{eq:17} uniformly in $p$.

We next establish the uniform estimate
\begin{equation}\label{eq:18}
 \norm{w}_{C^{2,\alpha}}\leq C\norm{L_pw}_{C^{0,\alpha}},
 \qquad w\perp\mathcal H_2,\quad p\in[-n-1,-n].
\end{equation}
Write $p=-n-\delta$, $0\leq\delta\leq1$.
By \eqref{eq:1}, the eigenvalue of $L_p$ on $\mathcal H_\ell$ is
\[
 2n+\delta-\ell(\ell+n-2).
\]
For $\ell=0,1$ these eigenvalues are at least $2n$ and $n+1$, respectively; for $\ell\geq3$ they are at most $-n-2$.
Expansion in spherical harmonics therefore gives
\[
 \norm{w}_{L^2}\leq\frac1{n+1}\norm{L_pw}_{L^2}
 \qquad(w\perp\mathcal H_2).
\]
The global Schauder estimate on the sphere, with its constant uniform in this bounded range of $p$, is
\[
 \norm{w}_{C^{2,\alpha}}
 \leq C\big(\norm{L_pw}_{C^{0,\alpha}}+\norm{w}_{C^0}\big).
\]
Using the interpolation inequality
$\norm{w}_{C^0}\leq\eta\norm{w}_{C^{2,\alpha}}+C_\eta\norm{w}_{L^2}$
and choosing $\eta$ small enough to absorb the first term proves \eqref{eq:18}.

Finally, $L_p$ commutes with $\Pi_2$.
Applying $I-\Pi_2$ to \eqref{eq:16}, and then using \eqref{eq:17}--\eqref{eq:18}, gives
\[
 \norm{(I-\Pi_2)v}_{C^{2,\alpha}}
 \leq C\norm{(I-\Pi_2)N_p(v)}_{C^{0,\alpha}}
 \leq C\norm{v}_{C^{2,\alpha}}^2.
\]
Together with \eqref{eq:15}, this yields
$\norm{v}_{C^{2,\alpha}}\leq C\norm{v}_{C^{2,\alpha}}^2$.
Choose $\rho>0$ smaller than both the radius used for the Taylor estimates and $C^{-1}$.
Then $\norm{v}_{C^{2,\alpha}}<\rho$ forces $v=0$.
The proof is complete.
\end{proof}

\section{Nonuniqueness}\label{sec:4}

Throughout this section, $n\geq8$.
We construct a local family of solutions near $h=1$ at
\[
 p_0=-2n-3.
\]
The kernel of the linearized equation consists of spherical harmonics of degree three.
We first impose simplex symmetry to obtain a one-dimensional kernel, and then complete the proof by determining the direction of the resulting branch.
Details are as follows.

We identify $\R^n$ isometrically with
\[
 H=\left\{y\in\R^{n+1}:\sum_{i=1}^{n+1}y_i=0\right\},
 \qquad S=H\cap\mathbb S^n.
\]
The vectors
\[
 \mathbf{v}_i=\sqrt{\frac{n+1}{n}}\left(e_i-\frac1{n+1}\one\right),
 \qquad 1\leq i\leq n+1,
\]
are the vertices of a regular $n$-simplex centered at the origin.
Here $e_i$ is the $i$th standard basis vector of $\R^{n+1}$ and $\one=(1,\ldots,1)^T$.
Its full symmetry group $G$ acts on $H$ by permutations of the coordinates.
Indeed, these permutations preserve the simplex, and every symmetry permutes its vertices, which span $H$.
We use the notation $\nabla$, $\Delta$, and $A_h$ on $S$ under the above isometry.
Integrals on $S$ use its round area measure, and $|S|$ denotes its area.

For $k\geq1$, put
\[
 s_k(y)=\sum_{i=1}^{n+1}y_i^k,
\]
and use the same notation for the restriction to $S$.
Thus $s_1=0$ and $s_2=1$ on $S$.
The orthogonal projection onto $H$ is
\[
 \Pi_H=I_{n+1}-\frac1{n+1}\one\otimes\one.
\]
Applying this projection to the Euclidean derivatives gives
\begin{equation}\label{eq:21}
 \begin{split}
 \Delta_Hs_k&=\frac{n}{n+1}k(k-1)s_{k-2},\qquad k\geq3,\\
 \nabla_Hs_k\cdot\nabla_Hs_l
 &=kl\left(s_{k+l-2}-\frac{s_{k-1}s_{l-1}}{n+1}\right),\qquad k,l\geq2.
 \end{split}
\end{equation}
In particular, $s_3$ is harmonic on $H$, with $s_3(1,1,-2,0,\ldots,0)=-6$.
Define
\begin{equation}\label{eq:22}
 \phi=\alpha_n s_3|_S,\qquad
 \alpha_n=\left(\int_Ss_3^2\right)^{-1/2}.
\end{equation}
Then $\int_S\phi^2=1$, and \eqref{eq:1} gives
\begin{equation}\label{eq:23}
 \Delta\phi=-3(n+1)\phi,\qquad
 \mathcal H_3^G=\SPAN\{\phi\}.
\end{equation}
To prove the second assertion, extend any $G$-invariant homogeneous harmonic cubic on $H$ to $\R^{n+1}$ by composing with $\Pi_H$.
This extension is a symmetric homogeneous cubic, hence a linear combination of $s_3$, $s_1s_2$, and $s_1^3$.
Its restriction to $H$ is therefore a multiple of $s_3$.

Fix $\alpha\in(0,1)$, and write
\[
 X=C^{2,\alpha}(S)^G,\qquad Y=C^{0,\alpha}(S)^G,
\]
where the superscript denotes $G$-invariant functions.
We shall use the closed subspaces
\[
 X_\perp=\left\{v\in X:\int_Sv\phi=0\right\},\qquad
 Y_\perp=\left\{g\in Y:\int_Sg\phi=0\right\}.
\]
The following two assertions give the required solutions.

\begin{lemma}\label{lem:6}
There are smooth functions $\beta(t)\in\R$ and $h(t)\in X$, defined for $|t|$ sufficiently small, such that
\begin{equation}\label{eq:24}
 \det A_{h(t)}=h(t)^{p_0+\beta(t)-1},\qquad
 \int_S(h(t)-1)\phi=t.
\end{equation}
Moreover, $\beta$ is even and
\begin{equation}\label{eq:25}
 h(t)=1+t\phi+t^2\varphi+O(t^3),\qquad
 \beta(t)=C_nt^2+O(t^4),
\end{equation}
where $\varphi\in X_\perp$ and $C_n\in\R$.
The first remainder is taken in $C^{2,\alpha}(S)^G$.
For every sufficiently small $t\ne0$, $h(t)$ is a nonconstant positive smooth support function with $A_{h(t)}>0$.
\end{lemma}

\begin{lemma}\label{lem:7}
The coefficient in \eqref{eq:25} is
\begin{equation}\label{eq:26}
 C_n=\frac{36\alpha_n^2(n+2)(2n^4+11n^3-110n^2-605n-570)}
 {(n+5)(n+6)(n+7)(n+8)(n+10)}.
\end{equation}
In particular, $C_n>0$ for $n\geq8$.
\end{lemma}

Assuming these for the moment, we prove the theorem.

\begin{proof}[Proof of Theorem~\ref{thm:2}]
Take the family in Lemma~\ref{lem:6}.
By \eqref{eq:25} and Lemma~\ref{lem:7}, $\beta(t_0)>0$ for every sufficiently small $t_0>0$.
Fix such a $t_0$ within the range of Lemma~\ref{lem:6}, and choose
\[
 0<\beta_n<\min\{\beta(t_0),n+3\}.
\]
Since $\beta$ is continuous and $\beta(0)=0$, the intermediate value theorem gives, for every $p\in(p_0,p_0+\beta_n)$, a $t\in(0,t_0)$ such that $\beta(t)=p-p_0$.
By \eqref{eq:24}, the corresponding $h(t)$ is a nonconstant positive solution of \eqref{LpMP1}.
It is invariant under $G$, the full symmetry group of the regular simplex described above.
This proves the theorem.
\end{proof}

It remains to prove Lemma~\ref{lem:6} and Lemma~\ref{lem:7}.
We first construct the local branch, and then compute its exponent to second order.

\subsection{A local reduction near the unit ball}

For $v$ in a sufficiently small neighborhood of zero in $X$, define
\begin{equation}\label{eq:27}
 F(p,v)=\det(I+A_v)-(1+v)^{p-1}.
\end{equation}
We choose this neighborhood so that $1+v>0$ and $I+A_v>0$.
The map $F$ takes values in $Y$ and is smooth in $(p,v)$, since the determinant is polynomial and the power is smooth on positive functions in the H\"older algebra.
Moreover, $F(p,0)=0$ and
\begin{equation}\label{eq:28}
 D_vF(p_0+\beta,0)=L-\beta I,\qquad L=\Delta+3(n+1).
\end{equation}
By \eqref{eq:1} and \eqref{eq:23}, the kernel of $L$ in $X$ is $\SPAN\{\phi\}$, and
\begin{equation}\label{eq:29}
 L:X_\perp\longrightarrow Y_\perp
\end{equation}
is an isomorphism.
Indeed, the eigenvalue on $\mathcal H_\ell$ is $3(n+1)-\ell(\ell+n-2)$.
For $\ell\ne3$, its absolute value is at least $n+3$.
Thus spherical harmonic expansion gives a unique $H^2$ solution orthogonal to $\mathcal H_3$ for every $g\in Y_\perp$.
Here $g$ is orthogonal to all of $\mathcal H_3$: its projection onto this space is $G$-invariant and hence a multiple of $\phi$.
Elliptic regularity gives a $C^{2,\alpha}$ solution, and uniqueness makes it $G$-invariant.
The spectral $L^2$ bound and the Schauder estimate, with interpolation as in \eqref{eq:18}, give
\[
 \norm{w}_{C^{2,\alpha}}\leq C\norm{Lw}_{C^{0,\alpha}},\qquad w\in X_\perp.
\]
This proves the boundedness of the inverse in \eqref{eq:29}.

\begin{proof}[Proof of Lemma~\ref{lem:6}]
For small $(t,\beta)$, consider the augmented equation
\begin{equation}\label{eq:30}
 F(p_0+\beta,t\phi+w)+\mu\phi=0,\qquad
 (w,\mu)\in X_\perp\times\R.
\end{equation}
Its derivative in $(w,\mu)$ at the origin is
\[
 (\widehat w,\widehat\mu)\longmapsto L\widehat w+\widehat\mu\phi.
\]
By \eqref{eq:29} and $\int_S\phi^2=1$, this is an isomorphism from $X_\perp\times\R$ to $Y$.
The implicit function theorem therefore gives unique smooth functions $w=w(t,\beta)\in X_\perp$ and $\mu=\mu(t,\beta)$ solving \eqref{eq:30} near the origin.
Since $F(p,0)=0$, uniqueness gives
\begin{equation}\label{eq:31}
 w(0,\beta)=0,\qquad \mu(0,\beta)=0.
\end{equation}
Differentiating \eqref{eq:30} in $t$ at $t=0$, we obtain
\[
 (L-\beta I)w_t(0,\beta)+(\mu_t(0,\beta)-\beta)\phi=0.
\]
The two terms lie in $Y_\perp$ and $\SPAN\{\phi\}$, respectively.
For small $\beta$, $L-\beta I$ remains invertible on the complementary spaces by \eqref{eq:29} and the Neumann series.
Consequently,
\begin{equation}\label{eq:32}
 w_t(0,\beta)=0,\qquad \mu_t(0,\beta)=\beta.
\end{equation}
In particular, $w(t,\beta)=O(t^2)$ uniformly for small $\beta$.

To obtain a solution of the original equation, we require $\mu=0$.
By \eqref{eq:31}, we can write
\[
 \mu(t,\beta)=t\Theta(t,\beta),\qquad
 \Theta(t,\beta)=\int_0^1\mu_t(\tau t,\beta)\dd \tau.
\]
The function $\Theta$ is smooth, and \eqref{eq:32} gives
\[
 \Theta(0,\beta)=\beta,\qquad
 \Theta(0,0)=0,\qquad \partial_\beta\Theta(0,0)=1.
\]
The scalar implicit function theorem yields a unique smooth $\beta=\beta(t)$ with $\beta(0)=0$ and $\Theta(t,\beta(t))=0$.
Set
\[
 h(t)=1+t\phi+w(t,\beta(t)).
\]
Then \eqref{eq:30} and $w(t,\beta)\in X_\perp$ prove both identities in \eqref{eq:24}.

We next determine the parity of $\beta$.
Define $R$ by $Rf(y)=f(-y)$.
The map $R$ preserves $X$, $Y$, and their complementary spaces, since antipodal reflection commutes with $G$ and $R\phi=-\phi$.
Moreover, $F(p,Rv)=RF(p,v)$.
Reflecting \eqref{eq:30} and using uniqueness gives
\[
 w(-t,\beta)=Rw(t,\beta),\qquad
 \mu(-t,\beta)=-\mu(t,\beta).
\]
Thus $\Theta$ is even in $t$, and uniqueness in the scalar implicit function theorem implies $\beta(-t)=\beta(t)$.
Taylor's formula and \eqref{eq:32} now give \eqref{eq:25}, with $\varphi=\frac12w_{tt}(0,0)\in X_\perp$.

For $t\ne0$, the second identity in \eqref{eq:24} rules out a constant $h(t)$, since $\phi$ is odd.
For small $t$, $h(t)>0$ and $A_{h(t)}>0$ by convergence to $1$ in $C^{2,\alpha}$.
Thus $h(t)$ is the support function of a strictly convex body containing the origin in its interior.
Elliptic regularity gives smoothness, as in the preliminaries.
The proof is complete.
\end{proof}

\subsection{The coefficient of the exponent}

We first calculate the polynomial moments needed below.
Write
\[
 \langle g\rangle=\frac1{|S|}\int_Sg.
\]

\begin{lemma}\label{lem:8}
The following identities hold:
\begin{equation}\label{eq:33}
 \begin{split}
 \langle s_4\rangle&=\frac{3n}{(n+1)(n+2)},\\
 \langle s_3^2\rangle&=\frac{6(n-1)}{(n+1)(n+2)(n+4)},\\
 \langle s_4^2\rangle&=\frac{3(3n^3+32n^2-8n+8)}{(n+1)^2(n+2)(n+4)(n+6)},\\
 \langle s_3s_5\rangle&=\frac{60n(n-1)}{(n+1)^2(n+2)(n+4)(n+6)},\\
 \langle s_3^2s_4\rangle&=\frac{18(n-1)(n^2+24n-12)}{(n+1)^2(n+2)(n+4)(n+6)(n+8)},\\
 \langle s_3^4\rangle&=\frac{108(n-1)(n^2+29n-42)}{(n+1)^2(n+2)(n+4)(n+6)(n+8)(n+10)}.
 \end{split}
\end{equation}
In particular,
\begin{equation}\label{eq:34}
 \alpha_n^2=\frac{(n+1)(n+2)(n+4)}{6(n-1)|S|}.
\end{equation}
\end{lemma}

\begin{proof}
For a homogeneous polynomial $P$ of degree $\ell$ on $H$, the polar-coordinate formula gives
\[
 (\Delta_HP)|_S=\Delta(P|_S)+\ell(\ell+n-2)P|_S.
\]
Since the spherical Laplacian has integral zero, we have
\begin{equation}\label{eq:35}
 \langle\Delta_HP\rangle=\ell(\ell+n-2)\langle P\rangle.
\end{equation}
For $P=s_4$ and $P=s_3^2$, respectively, \eqref{eq:21} and \eqref{eq:35} give
\[
 \begin{split}
 4(n+2)\langle s_4\rangle&=\frac{12n}{n+1},\\
 6(n+4)\langle s_3^2\rangle&=18\left(\langle s_4\rangle-\frac1{n+1}\right).
 \end{split}
\]
These prove the first two identities in \eqref{eq:33}.
The same calculation with $P=s_6$ gives
\[
 \langle s_6\rangle=\frac{15n^2}{(n+1)^2(n+2)(n+4)}.
\]
To obtain the remaining moments, apply the product rule in \eqref{eq:21} and then \eqref{eq:35} to $s_4^2$, $s_3s_5$, $s_3^2s_4$, and $s_3^4$.
This yields
\[
 \begin{split}
 8(n+6)\langle s_4^2\rangle
 &=\frac{24n}{n+1}\langle s_4\rangle+32\langle s_6\rangle
       -\frac{32}{n+1}\langle s_3^2\rangle,\\
 8(n+6)\langle s_3s_5\rangle
 &=\frac{20n}{n+1}\langle s_3^2\rangle+30\langle s_6\rangle
       -\frac{30}{n+1}\langle s_4\rangle,\\
 10(n+8)\langle s_3^2s_4\rangle
 &=18\langle s_4^2\rangle-\frac{18}{n+1}\langle s_4\rangle
       +\frac{12(n-4)}{n+1}\langle s_3^2\rangle+48\langle s_3s_5\rangle,\\
 12(n+10)\langle s_3^4\rangle
 &=108\left(\langle s_3^2s_4\rangle-\frac1{n+1}\langle s_3^2\rangle\right).
 \end{split}
\]
Substitution in this order proves \eqref{eq:33}.
Finally, \eqref{eq:34} follows from \eqref{eq:22} and the value of $\langle s_3^2\rangle$.
\end{proof}

We now use these moments to evaluate the scalar equation along the local branch.

\begin{proof}[Proof of Lemma~\ref{lem:7}]
All derivatives of $F$ in this proof are evaluated at $(p_0,0)$.
Write
\[
 Q(f,g)=(D_v^2F)[f,g],\qquad T(f,g,u)=(D_v^3F)[f,g,u].
\]
Since $p_0-1=-2(n+2)$, differentiation of \eqref{eq:27} gives
\begin{equation}\label{eq:36}
 Q(f,g)=(\TR A_f)(\TR A_g)-\TR(A_fA_g)-2(n+2)(2n+5)fg.
\end{equation}

We first calculate the matrices generated by $s_3$.
For $y\in S$, let
\[
 P_y=\Pi_H-yy^T,\qquad
 \Lambda_y=\DIAG(y_1,\ldots,y_{n+1}),\qquad
 B=(P_y\Lambda_yP_y)|_{T_yS}.
\]
Here $P_y$ is the orthogonal projection onto $T_yS=H\cap y^\perp$.
The identity matrix and traces below are on this $(n-1)$-dimensional tangent space.
Since $s_3$ is homogeneous of degree three and its ambient Hessian is $6\Lambda_y$, we have
\[
 A_{s_3}=6B-2s_3I.
\]
To compute the traces of $B$, put $\mathbf{e}=(n+1)^{-1/2}\one$ and $U=\mathbf{e}\mathbf{e}^T+yy^T$, so that $P_y=I_{n+1}-U$.
Using $s_1=0$, $s_2=1$, and cyclicity of the ambient trace gives
\[
 \begin{split}
 \TR B&=-s_3,\\
 \TR(B^2)&=1-2\left(\frac1{n+1}+s_4\right)
                 +\frac2{n+1}+s_3^2
             =1-2s_4+s_3^2,\\
 \TR(B^3)&=s_3-3\left(\frac{s_3}{n+1}+s_5\right)
          +3\left(\frac{2s_3}{n+1}+s_3s_4\right)
          -\left(\frac{3s_3}{n+1}+s_3^3\right)\\
          &=s_3-3s_5+3s_3s_4-s_3^3.
 \end{split}
\]
For example, the last three terms in the cubic calculation are
$-3\TR(\Lambda_y^3U)$, $3\TR(\Lambda_y^2U\Lambda_yU)$, and $-\TR((\Lambda_yU)^3)$.
Expanding the powers of $6B-2s_3I$ therefore yields
\begin{equation}\label{eq:37}
 \begin{split}
 \TR A_{s_3}&=-2(n+2)s_3,\\
 \TR(A_{s_3}^2)&=36-72s_4+4(n+14)s_3^2,\\
 \TR(A_{s_3}^3)&=-648s_5+1080s_3s_4-8(n+62)s_3^3.
 \end{split}
\end{equation}

The coefficient of $t^2$ in \eqref{eq:24} is
\[
 L\varphi+\frac12Q(\phi,\phi)=0.
\]
By \eqref{eq:36}--\eqref{eq:37}, this becomes
\begin{equation}\label{eq:38}
 L\varphi=\alpha_n^2\left(18-36s_4+3(n+10)s_3^2\right).
\end{equation}
The span of $1,s_4,s_3^2$ is preserved by $L$.
Indeed, \eqref{eq:21} and homogeneity give
\[
 \begin{split}
 L1&=3(n+1),\\
 Ls_4&=\frac{12n}{n+1}-(n+5)s_4,\\
 L(s_3^2)&=18s_4-\frac{18}{n+1}-3(n+7)s_3^2.
 \end{split}
\]
Matching coefficients in \eqref{eq:38}, we obtain
\begin{equation}\label{eq:39}
 \begin{split}
 \varphi&=\alpha_n^2\widetilde\varphi,\\
 \widetilde\varphi
 &=\frac{6(n^2-n-15)}{(n+1)(n+5)(n+7)}
   +\frac{18(n+4)}{(n+5)(n+7)}s_4-\frac{n+10}{n+7}s_3^2.
 \end{split}
\end{equation}
The right-hand side is even and hence orthogonal to $\phi$.
It is therefore the unique correction in $X_\perp$ by \eqref{eq:29}.

We next extract the scalar condition at order three.
Taylor's formula and the second equality in \eqref{eq:24} give
\[
 h(t)=1+t\phi+t^2\varphi+t^3\eta+O(t^4),\qquad \eta\in X_\perp.
\]
To keep track of the varying exponent, Taylor expansion of \eqref{eq:27} gives
\[
 F(p_0+\beta,v)=Lv-\beta v+\frac12Q(v,v)+\frac16T(v,v,v)
 +O\left(\norm{v}_{C^{2,\alpha}}^4+|\beta|\norm{v}_{C^{2,\alpha}}^2
                  +\beta^2\norm{v}_{C^{2,\alpha}}\right)
\]
in $C^{0,\alpha}$ near $(\beta,v)=(0,0)$.
Since $\beta(t)=C_nt^2+O(t^4)$, the coefficient of $t^3$ in the equation is
\[
 L\eta-C_n\phi+Q(\phi,\varphi)+\frac16T(\phi,\phi,\phi)=0.
\]
Pairing with $\phi$, using self-adjointness of $L$, $L\phi=0$, and $\int_S\phi^2=1$, gives
\begin{equation}\label{eq:40}
 C_n=\int_S\phi Q(\phi,\varphi)+\frac16\int_S\phi T(\phi,\phi,\phi).
\end{equation}

An integration-by-parts identity simplifies the mixed term.
For every smooth $g$, the Codazzi identity for $A_g$ gives
\[
 \nabla_i\big((\TR A_g)\delta_{ij}-(A_g)_{ij}\big)=0.
\]
Consequently,
\[
 \begin{split}
 &\int_S f\big((\TR A_g)(\TR A_u)-\TR(A_gA_u)\big)\\
 &\quad=\int_S f\big((\TR A_g)\delta_{ij}-(A_g)_{ij}\big)(u_{,ij}+u\delta_{ij})\\
 &\quad=-\int_S\big((\TR A_g)\delta_{ij}-(A_g)_{ij}\big)f_{,i}u_{,j}
           +(n-2)\int_Sfu\TR A_g.
 \end{split}
\]
The last expression is symmetric in $f,u$, while the first is symmetric in $g,u$.
It is therefore symmetric in all three functions.
The last term of \eqref{eq:36} has the same property after multiplication by $f$ and integration.
Using $Q(\phi,\phi)=-2L\varphi$, we obtain
\begin{equation}\label{eq:41}
 \int_S\phi Q(\phi,\varphi)=\int_S\varphi Q(\phi,\phi)
 =-2\int_S\varphi L\varphi.
\end{equation}

For the pure cubic term, the third derivative of the determinant and the power gives
\[
 \begin{split}
 T(s_3,s_3,s_3)
 &=(\TR A_{s_3})^3-3(\TR A_{s_3})\TR(A_{s_3}^2)+2\TR(A_{s_3}^3)\\
 &\quad+2(n+2)(2n+5)(2n+6)s_3^3.
 \end{split}
\]
By \eqref{eq:37},
\begin{equation}\label{eq:42}
 \begin{split}
 \gamma:=\frac16\langle s_3T(s_3,s_3,s_3)\rangle
 &=36(n+2)\langle s_3^2\rangle-72(n-3)\langle s_3^2s_4\rangle\\
 &\quad-216\langle s_3s_5\rangle+(6n^2+70n-44)\langle s_3^4\rangle.
 \end{split}
\end{equation}
Combining \eqref{eq:40}--\eqref{eq:42} with $\varphi=\alpha_n^2\widetilde\varphi$ and $\phi=\alpha_ns_3$, we have
\begin{equation}\label{eq:43}
 C_n=\alpha_n^4|S|\left(-2\langle\widetilde\varphi L\widetilde\varphi\rangle+\gamma\right)
 =\frac{\alpha_n^2}{\langle s_3^2\rangle}
       \left(-2\langle\widetilde\varphi L\widetilde\varphi\rangle+\gamma\right).
\end{equation}
Here we used $\alpha_n^2|S|\langle s_3^2\rangle=1$, so the coefficient retains the $L^2(\dd x)$ normalization in \eqref{eq:22}.

It remains to substitute the moments in Lemma~\ref{lem:8}.
For completeness, write $\widetilde\varphi=a+bs_4+cs_3^2$, with the coefficients in \eqref{eq:39}.
By \eqref{eq:38},
\[
 \begin{split}
 \langle\widetilde\varphi L\widetilde\varphi\rangle
 &=a\left(18-36\langle s_4\rangle+3(n+10)\langle s_3^2\rangle\right)\\
 &\quad+b\left(18\langle s_4\rangle-36\langle s_4^2\rangle
                         +3(n+10)\langle s_3^2s_4\rangle\right)\\
 &\quad+c\left(18\langle s_3^2\rangle-36\langle s_3^2s_4\rangle
                         +3(n+10)\langle s_3^4\rangle\right)\\
 &=\frac{108(n-1)(n^4+21n^3+139n^2+193n-426)}
 {(n+1)(n+4)(n+5)(n+6)(n+7)(n+8)}.
 \end{split}
\]
Likewise, \eqref{eq:42} becomes
\[
 \gamma=\frac{216(n-1)(n^3+21n^2+73n-138)}
 {(n+1)(n+4)(n+6)(n+8)(n+10)}.
\]
Substitution in \eqref{eq:43} proves \eqref{eq:26}.
For its sign, write $n=8+r$, where $r\geq0$.
Then
\[
 2n^4+11n^3-110n^2-605n-570
 =2r^4+75r^3+922r^2+3843r+1374>0.
\]
All remaining factors in \eqref{eq:26} are positive.
Thus $C_n>0$ for $n\geq8$, as asserted.
\end{proof}

\noindent\textbf{Declaration of AI Use.}
During the preparation, the authors used OpenAI's ChatGPT-5.6 to assist with language refinement, mathematical exposition, and the examination of selected proofs and calculations. All AI-assisted content was critically reviewed and independently verified by the authors, who take full responsibility for the accuracy and integrity of the manuscript.


\end{document}